\documentclass{article}

\newtheorem{thm}{Theorem}[section]

\newtheorem{lem}[thm]{Lemma}

\begin{document}

\title{{\bf An Improved Error Bound for Multiquadric Interpolation}}         
\author{{\bf Lin-Tian Luh} \\Department of Mathematics, Providence University\\ Shalu, Taichung\\ email:ltluh@pu.edu.tw \\ phone:(04)26328001 ext. 15126 \\ fax:(04)26324653}        
\maketitle
{\bf Abstract.}  For multiquadric interpolation it's well known that the most powerful error bound is the so-called exponential-type error bound. It's of the form $|f(x)-s(x)|\leq C\lambda ^{1/d}\| f\| _{h}$ where $0<\lambda <1$ is a constant, C is a constant, d is the fill distance which roughly speaking measures the spacing of the data points, $s(x)$ is the interpolating function of $f(x)$, and $\| f\| _{h}$ denotes the norm of f induced by the multiquadric or inverse multiquadric. This error bound for $x \in R^{n}$ was put forward by Madych and Nelson in 1992 and converges to zero very fast as $d\rightarrow 0$. The drawback is that both C and $\lambda$ get worse as $n\rightarrow \infty $. In particular, $\lambda \rightarrow 1^{-}$ very fast as $n\rightarrow \infty$. In this paper, we raise an error bound of the form $\ |f(x)-s(x)| \leq C'\sqrt{d} \lambda '^{1/d}\| f\| _{h} $ where $\lambda'$ can be independent of the dimension n and $0<\lambda '< \lambda < 1 $. Moreover, $C'$ is only slightly different from C. What's noteworthy is that both $C'$ and $\lambda '$ can be computed without slight difficulty. 
\section{Introduction}       
We begin with some fundamental knowledge about radial basis functions. Let $h$ be a continuous function in $R^{n}$ and be strictly conditionally positive definite of order $m$. Given data points $(x_{j}, f_{j}),j=1,\ldots ,N$, where $X=\{ x_{1},\ldots ,x_{N}\} $ is a subset of points in $R^{n}$ and the $f_{j}'$s are real or complex numbers, the so-called $h$ spline interpolant of these data is the function $s$ defined by
\begin{equation}
  s(x)=p(x)+\sum_{j=1}^{N}c_{j}h(x-x_{j}),
\end{equation} 
where $p(x)$ is a polynomial in $P_{m-1}^{n}$, the space of n-variable polynomials of degree not more than $m-1$, and the $c_{j}$'s are chosen so that 
\begin{equation}
  \sum_{j=1}^{N}c_{j}q(x_{j})=0
\end{equation}
for all polynomials $q$ in $P_{m-1}^{n}$ and
\begin{equation}
  p(x_{i})+\sum_{j=1}^{N}c_{j}h(x_{i}-x_{j})=f_{i},\ \ i=1,\ldots , N.
\end{equation}
It is well known that the linear system (2) and (3) has a unique solution if $X$ is a determining set for $P_{m-1}^{n}$ and $h$ is strictly conditionally positive definite. Further details can be found in \cite{MN1}. Therefore in our case the interpolant $s(x)$ is well defined.

Here, $X$ is a determining set for $P_{m-1}^{n}$ means that $X$ does not lie in the zero set of any nontrivial member of $P_{m-1}^{n}$.

In this paper $h$ is defined by
\begin{equation}
  h(x):=\Gamma (-\frac{\beta}{2})(c^{2}+|x|^{2})^{\frac{\beta}{2}},\ \ \beta\in R\backslash 2N_{\geq 0},\ \ c>0.
\end{equation}
where $|x|$ is the Euclidean norm of $x$, $\Gamma$ is the classical gamma function and $\beta,\ c$ are constants. The function $h$ is called multiquadric or inverse multiquadric, respectively, depending on $\beta >0$, or $\beta<0$.

In \cite{MN3} Madych and Nelson raise the famous exponential-type error bound for scattered data interpolation of multiquadric and inverse multiquadric, as mentioned in the abstract. The computation of C and the crucial constant $\lambda$ is solved by Luh in \cite{Lu3}. This brings the error bound forward to the possibility of application. However, in practical application, too many data points may result in a big condition number for solving the linear system (2) and (3). A good solution is to improve the error bound so that a satisfactory error estimate can be reached before ill-conditioning occurs. This is the purpose of this paper.
\subsection{Some Fundamental Knowledge about Polynomials}    
Let $P_{l}^{n}$ be the space of polynomials of degree not exceeding $l$ in n variables and $E\subseteq R^{n}$ be compact. It's well known that $dimP_{l}^{n}=\left( \begin{array}{c}
                                                         n+l \\
                                                          n
                                                      \end{array} \right) $. In this section we will denote $dimP_{l}^{n}$ by $N$. It's also well known that if $x_{1},\ldots ,x_{N} \in E$ do not lie on the zero set of any nontrivial $q\in P_{l}^{n}$, there exist Lagrange polynomials $l_{i}$ of degree $l$ defined by $l_{i}(x_{j})=\delta_{ij},\ 1\leq i, j \leq N$. This guarantees that for any $f\in C(E),\ (\Pi _{l}f)(x):=\sum _{i=1}^{N}f(x_{i})l_{i}(x)$ is its interpolating polynomial. It's easily seen that $\Pi _{l}(p)=p$ for all $p\in P_{l}^{n}$ and hence the mapping $\Pi _{l}:C(E)\longrightarrow P_{l}^{n}$ is a projection. Let 
$$\| \Pi _{l}\| := \max_{x\in E}\sum _{i=1}^{N}|l_{i}(x)|.$$
It's easy to show that for any $p\in P_{l}^{n}$,
$$\| p\| _{\infty}\leq \| \Pi _{l}\| \max_{1\leq i\leq N}|p(x_{i})|$$
if the domain of $p$ is $E$.
In this paper we will discuss mainly the case $E=T_{n}$, an n-dimensional simplex whose definition can be found in \cite{Fl}.

Now we are going to adopt {\bf barycentric coordinates}. Let $v_{i},\ 1\leq i\leq n+1$, be the vertices of $T_{n}$. Then any $x\in T_{n}$ can be written as a convex combination of the vertices:
$$x=\sum_{i=1}^{n+1}c_{i}v_{i},\ \sum_{i=1}^{n+1}c_{i}=1,\ c_{i}\geq 0.$$
We define ``{\bf equally spaced}'' points of {\bf degree} $l$ to be those points whose barycentric coordinates are of the form
$$(k_{1}/l,k_{2}/l,\ldots ,k_{n+1}/l),\ k_{i}\ nonnegative\  integers \ with \ \sum_{i=1}^{n+1}k_{i}=l.$$
It's obvious that the number of such points is exactly $N=dimP_{l}^{n}$. Moreover, by \cite{Bo}, we know that equally spaced points form a determining set for $P_{l}^{n}$. We will need the following lemma which we cite from \cite{Bo}.

\begin{lem}
  For the above equally spaced points $\| \Pi_{l}\| \leq \left( \begin{array}{c}
                                                                   2l-1\\ l
                                                                \end{array}  \right) $. Moreover, as $n\rightarrow \infty,\ \| \Pi_{l}\| \longrightarrow \left( \begin{array}{c}
                                                                           2l-1\\ l
                                                                        \end{array}  \right)$. 
\end{lem}
Now, we are going to prove a lemma which plays a crucial role in our construction of the error bound.
\begin{lem}
  Let $Q\subseteq R^{n}$ be a simplex and $Y$ be the set of equally spaced points of degree $l$ in $Q$. Then, for any point $x$ in $Q$, there is a measure $\sigma$ supported on $Y$ such that 
$$\int p(y)d\sigma (y)=p(x)$$
for all $p$ in $P_{l}^{n}$, and 
$$\int d|\sigma |(y)\leq  \left( \begin{array}{c}
                            2l-1\\ l
                         \end{array} \right) $$
\end{lem}
{\bf Proof}. Let $Y=\{ y_{1},\ldots ,y_{N}\} $ be the set of equally spaced points of degree $l$ in $Q$. Denote $P_{l}^{n}$ by $V$. For any $x\in Q$, let $\delta_{x}$ be the point-evaluation functional. Define $T:V\longrightarrow T(V)\subseteq R^{N}$ by $T(v)=(\delta_{y_{i}}(v))_{y_{i}\in Y}$. Then $T$ is injective. Define $\tilde{\psi}$ on $T(V)$ by $\tilde{\psi}(w)=\delta _{x}(T^{-1}w)$. By the Hahn-Banach theorem, $\tilde{\psi}$ has a norm-preserving extension $\tilde{\psi}_{ext}$ to $R^{N}$. By the Riesz representation theorem, each linear functional on $R^{N}$ can be represented by the inner product with a fixed vector. Thus, there exists $z\in R^{N}$ with
$$\tilde{\psi}_{ext}(w)=\sum _{j=1}^{N}z_{j}w_{j}$$
and $\| z\| _{(R^{N})^{*}}= \| \tilde{\psi}_{ext}\| $. If we adopt the $l_{\infty}$-norm on $R^{N}$, the dual norm will be the $l_{1}$-norm. Thus $\| z\| _{(R^{N})^{*}}=\| z\| _{1}= \| \tilde{\psi}_{ext}\| =\| \tilde{\psi} \| =\| \delta _{x}T^{-1}\|$.

Now, for any $p\in V$, by setting $w=T(p)$, we have 
$$\delta_{x}(p)=\delta _{x}(T^{-1}w)=\tilde{\psi}(w)=\tilde{\psi}_{ext}(w)=\sum _{j=1}^{N}z_{j}w_{j}=\sum _{j=1}^{N}z_{j}\delta_{y_{j}}(p).$$
This gives
\begin{equation}
  p(x)=\sum _{j=1}^{N}z_{j}p(y_{j})
\end{equation}
where $|z_{1}|+\cdots +|z_{N}|= \| \delta_{x}T^{-1}\|$.

Note that
\begin{eqnarray*}
  \|\delta_{x}T^{-1}\| & =  & \sup _{\begin{array}{c}
                                       w\in T(V)\\ w \neq 0
                                     \end{array} } \frac{\| \delta_{x}T^{-1}(w)\|}{\| w\| _{R^{N}}} \\
                       & =  & \sup _{\begin{array}{c}
                                       w\in T(V)\\ w \neq 0
                                     \end{array}} \frac{|\delta_{x}p|}{\| T(p)\| _{R^{N}}}\\
                       & \leq & \sup _{\begin{array}{c}
                                         p\in V\\ p\neq 0
                                       \end{array}} \frac{|p(x)|}{\max _{j=1,\ldots ,N}|p(y_{j})|}\\
                       & \leq & \sup _{\begin{array}{c}
                                         p\in V\\ p\neq 0
                                       \end{array}} \frac{\| \Pi _{l}\| \max _{j=1,\ldots ,N}|p(y_{j})|}{\max _{j=1,\ldots ,N}|p(y_{j})|} \\
                       & = & \| \Pi _{l}\| \\
                       & \leq & \left( \begin{array}{c}
                                         2l-1 \\l   
                                       \end{array} \right).             
\end{eqnarray*}
Therefore $$|z_{1}|+\cdots +|z_{N}| \leq \left( \begin{array}{c}
                                                2l-1\\l
                                              \end{array} \right)$$ and our lemma follows immediately by (5) and letting $\sigma(\{ y_{j}\} )=z_{j},\ j=1,\ldots ,N$. \hspace{2cm} \ \ \ \ $\sharp$

\subsection{Some Fundamental Knowledge for Interpolation}
In our development of the error estimate some basic ingredients should be introduced. We emphasize on two objects:first, the function space where the interpolating and interpolated functions come from;second, an important lemma introduced by Gelfand and Vilenkin in \cite{GV}, but modified by Madych and Nelson in \cite{MN2}, which reflects deeply the features of a radial basis function.

Let $\cal D$ denote the space of complex-valued functions on $R^{n}$ that are compactly supported and infinitely differentiable. The Fourier transform of a function $\phi$ in $\cal D$ is $$\hat{\phi}(\xi) = \int e^{-i<x,\xi>}\phi(x)dx.$$
A continuous function $h$ is called {\bf conditionally positive definite} of order $m$ if 
$$\int h(x)\phi(x)\ast \tilde{\phi}(x)dx\geq 0$$ 
holds whenever $\phi=p(D)\psi$ with $\psi$ in $\cal D$ and $p(D)$ a linear homogeneous constant coefficient differential operator of order $m$. Here $\tilde{\phi}=\overline{\phi(-x)}$ and $\ast $ denotes the convolution product
$$\phi_{1}\ast \phi_{2}(t)=\int \phi_{1}(x)\phi_{2}(t-x)dx.$$
This definition, introduced in \cite{MN2}, looks complicated. However it's equivalent to the generally used definition. For further details we refer reader to \cite{MN1} and \cite{MN2}.

Now we have come to the afore-mentioned crucial lemma. It's as follows. If $h$ is a continuous conditionally positive definite function of order $m$, the Fourier transform of $h$ uniquely determines a positive Borel measure $\mu$ on $R^{n}\backslash \{ 0\} $ and constants $a_{\gamma},|\gamma|=2m$ as follows: For all $\psi\in \cal D$
\begin{eqnarray}
  \int h(x)\psi(x)dx & = & \int \left\{ \hat{\psi}(\xi)-\hat{\chi}(\xi)\sum_{|\gamma|<2m}D^{r}\hat{\psi}(0)\frac{\xi^{r}}{r!}\right\} d\mu (\xi)\\
                     &   & +\sum _{|\gamma|\leq 2m}D^{r}\hat{\psi}(0)\frac{a_{\gamma}}{r!}, \nonumber
\end{eqnarray}
where for every choice of complex numbers $c_{\alpha},|\alpha|=m$,
$$\sum _{|\alpha|=m}\sum _{|\beta|=m}a_{\alpha+\beta}c_{\alpha}\overline{c_{\beta}}\geq 0.$$
Here $\chi$ is a function in $\cal D$ such that $1-\hat{\chi}(\xi)$ has a zero of order $2m+1$ at $\xi =0$; both of the integrals $\int _{0<|\xi|<1}|\xi|^{2m}d\mu (\xi),\ \int _{|\xi|\geq 1} d\mu (\xi)$ are finite. The choice of $\chi$ affects the value of the coefficients $a_{\gamma}$ for $|\gamma|<2m$.

Our error estimate is based on a function space called native space whose characterizations can be found in \cite{Lu1}, \cite{Lu2} and \cite{We}. Although there are simpler expressions for the function space, we still adopt the one used by Madych and Nelson in \cite{MN2} to show the author's respect for them. First, let
$${\cal D}_{m}=\left\{ \phi \in {\cal D} : \int x^{\alpha }\phi (x)dx=0\ \ for\ all\ |\alpha |<m\right\},$$
Then the {\bf native space ${\cal C}_{h,m}$} is the class of those continuous functions $f$ which satisfy 
\begin{equation}
  \left| \int f(x)\phi(x)dx\right| \leq c(f)\left\{ \int h(x-y)\phi(x)\overline{\phi(y)}dxdy\right\} ^{1/2}
\end{equation}
for some constant $c(f)$ and all $\phi$ in ${\cal D}_{m}$. If $f\in {\cal C}_{h,m}$, let $\| f\| _{h}$ denote the smallest constant $c(f)$ for which (7) is true. Here $\| f\| _{h}$ is a semi-norm and ${\cal C}_{h,m}$ is a semi-Hilbert space; in the case $m=0$ it is a norm and a Hilbert space respectively.
\section{Main Results}
It's well known that the function $h$ defined in (4) is conditionally positive definite of order $m=0$ if $\beta<0$, and $m=\lceil \frac{\beta}{2}\rceil$ if $\beta>0$. With this in mind we have the following lemma.
\begin{lem}
  Let $h$ be as in (4) and $m$ be its order of conditional positive definiteness. There exists a positive constant $\rho$ such that 
\begin{equation}
  \int _{R^{n}} |\xi |^{l}d\mu (\xi)\leq (\sqrt{2})^{n+\beta +1}\cdot (\sqrt{\pi})^{n+1}\cdot n\alpha_{n}\cdot c^{\beta-l}\cdot \Delta _{0}\cdot \rho ^{l}\cdot l!
\end{equation}
for all integer $l\geq 2m+2$ where $\mu$ is defined in (6), $\alpha_{n}$ denotes the volume of the unit ball in $R^{n}$, $c$ is as in (4), and $\Delta _{0}$ is a positive constant.
\end{lem}
{\bf Proof}. Let ${\cal K}_{\nu}$ denote a modified Bessel function of the second kind. Then
\begin{eqnarray*}
  &   & \int _{R^{n}} |\xi |^{l}d\mu (\xi) \\
  & = & \int _{R^{n}}|\xi|^{l}2\pi ^{n/2}\left( \frac{|\xi|}{2c}\right)^{-\frac{n+\beta}{2}}{\cal K}_{\frac{n+\beta}{2}}(c|\xi|)d\xi  \\
  & = & 2\pi^{n/2}\left( \frac{1}{2c}\right) ^{-\frac{n+\beta}{2}}\int_{R^{n}}|\xi|^{l-\frac{n+\beta}{2}}{\cal K}_{\frac{n+\beta}{2}}(c|\xi|)d\xi \\
  & \sim & \sqrt{\frac{\pi}{2}}2\pi^{n/2}\left( \frac{1}{2c}\right) ^{-\frac{n+\beta}{2}}\int_{R^{n}}|\xi|^{l-\frac{n+\beta}{2}}\frac{1}{\sqrt{c|\xi|}e^{c|\xi|}}d\xi \\
  & = & \sqrt{\frac{\pi}{2}}2\pi^{n/2}\left( \frac{1}{2c}\right) ^{-\frac{n+\beta}{2}}n\alpha_{n}\int_{0}^{\infty}r^{l-\frac{n+\beta}{2}}\frac{r^{n-1}}{\sqrt{c|r|}e^{c|r|}}dr \\
  & = & \sqrt{\frac{\pi}{2}}2\pi^{n/2}(2c)^{\frac{n+\beta}{2}}\frac{1}{\sqrt{c}}n\alpha_{n}\int_{0}^{\infty}\frac{r^{l+\frac{n-\beta -3}{2}}}{e^{cr}}dr \\
  & = & \sqrt{\frac{\pi}{2}}2\pi^{n/2}(2c)^{\frac{n+\beta}{2}}\frac{1}{\sqrt{c}}n\alpha_{n}\frac{1}{c^{l+\frac{n-\beta-1}{2}}}\int_{0}^{\infty}\frac{r^{l+\frac{n-\beta-3}{2}}}{e^{r}}dr \\
  & = & 2^{\frac{n+\beta+1}{2}}\pi^{\frac{n+1}{2}}n\alpha_{n}c^{\beta-l}\int_{0}^{\infty}\frac{r^{l'}}{e^{r}}dr\ where\ l'=l+\frac{n-\beta-3}{2}.
\end{eqnarray*}
Note that if $\beta<0$, then $m=0$ and $l\geq 2m+2=2$. This implies $l'>0$. If $\beta>0$, then $m=\lceil \frac{\beta}{2}\rceil $ and $l\geq 2m+2=2\lceil \frac{\beta}{2} \rceil +2$. This implies $l'>0$. In any case $l'>0$.

Our proof then proceeds in three cases. Let $l''=\lceil l' \rceil $ which is the smallest integer greater than or equal to $l'$.

{\bf Case1.} Assume $l''>l$. Let $l''=l+s$. Then
$$\int_{0}^{\infty}\frac{r^{l'}}{e^{r}}dr\leq \int_{0}^{\infty}\frac{r^{l''}}{e^{r}}dr=l''!=(l+s)(l+s-1)\cdots (l+1)l!$$
and
$$\int_{0}^{\infty}\frac{r^{l'+1}}{e^{r}}dr\leq \int_{0}^{\infty}\frac{r^{l''+1}}{e^{r}}dr=(l''+1)!=(l+s+1)(l+s)\cdots(l+2)(l+1)l!.$$
Note that
$$\frac{(l+s+1)(l+s)\cdots (l+2)}{(l+s)(l+s-1)\cdots (l+1)}=\frac{l+s+1}{l+1}.$$
(i)Assume $\beta<0$. Then $m=0$ and $l\geq 2$. This gives
$$\frac{l+s+1}{l+1}\leq \frac{3+s}{3}.$$
Let $\rho =\frac{3+s}{3}$. Then 
$$\int_{0}^{\infty}\frac{r^{l''+1}}{e^{r}}dr\leq \Delta_{0}\rho^{l+1}(l+1)!$$
if $\int_{0}^{\infty}\frac{r^{l''}}{e^{r}}dr\leq \Delta_{0}\rho^{l}l!$. The smallest $l''$ is $l_{0}''=2+s$. Now,
\begin{eqnarray*}
\int_{0}^{\infty}\frac{r^{l_{0}''}}{e^{r}}dr & = & l_{0}''!=(2+s)(2+s-1)\cdots (3)l!\ \  where\ \  l=2 \\
                                             & = & \frac{(2+s)(2+s-1)\cdots (3)}{\rho^{2}}\cdot \rho^{l}l!(l=2) \\
                                             & = & \Delta_{0}\rho^{2}2!\ \ where\ \ \Delta_{0}=\frac{(2+s)(2+s-1)\cdots(3)}{\rho^{2}}. 
\end{eqnarray*}
It follows that $\int_{0}^{\infty}\frac{r^{l'}}{e^{r}}dr\leq \Delta_{0}\rho^{l}l!$ for all $l\geq 2$.

(ii)Assume $\beta>0$. Then $m=\lceil \frac{\beta}{2}\rceil $ and $l\geq 2m+2$. This gives
$$\frac{l+s+1}{l+1} \leq \frac{2m+3+s}{2m+3}=1+\frac{s}{2m+3}.$$
Let $\rho=1+\frac{s}{2m+3}$. Then
$$\int_{0}^{\infty}\frac{r^{l''+1}}{e^{r}}dr\leq \Delta_{0}\rho^{l+1}(l+1)!$$
if $\int_{0}^{\infty}\frac{r^{l''}}{e^{r}}dr\leq \Delta_{0}\rho^{l}l!$. The smallest $l''$ is $l_{0}''=2m+2+s$ when $l=2m+2$. Now,
\begin{eqnarray*}
 \lefteqn{\int_{0}^{\infty}\frac{r^{l_{0}''}}{e^{r}}dr}\\
   & = & l_{0}''!=(2m+2+s)(2m+1+s)\cdots(2m+3)(2m+2)!\\
   & = & \frac{(2m+2+s)(2m+1+s)\cdots (2m+3)}{\rho^{2m+2}}\cdot \rho^{2m+2}(2m+2)!\\
   & = & \Delta_{0}\rho^{2m+2}(2m+2)!\  where\ \Delta_{0}=\frac{(2m+2+s)(2m+1+s)\cdots (2m+3)}{\rho^{2m+2}}. 
\end{eqnarray*}
It follows that $\int_{0}^{\infty}\frac{r^{l'}}{e^{r}}dr \leq \Delta_{0}\rho^{l}l!$ for all $l\geq 2m+2$.

{\bf Case2.} Assume $l''<l$. Let $l''=l-s$ where $s>0$. Then
$$\int_{0}^{\infty}\frac{r^{l'}}{e^{r}}dr\leq \int_{0}^{\infty}\frac{r^{l''}}{e^{r}}dr=l''!=(l-s)!=\frac{1}{l(l-1)\cdots (l-s+1)}\cdot l!$$
and
\begin{eqnarray*}
  \int_{0}^{\infty}\frac{r^{l'+1}}{e^{r}}dr & \leq & \int_{0}^{\infty}\frac{r^{l''+1}}{e^{r}}dr\\
                                            & = & (l''+1)!=(l-s+1)!=\frac{1}{(l+1)l\cdots (l-s+2)}\cdot (l+1)!.
\end{eqnarray*}
Note that
\begin{eqnarray*}
  &   & \left\{ \frac{1}{(l+1)l\cdots (l-s+2)} /\frac{1}{l(l-1)\cdots (l-s+1)}\right\} \\
  & = & \frac{l(l-1)\cdots (l-s+1)}{(l+1)l\cdots (l-s+2)} \\
  & = & \frac{l-s+1}{l+1}
\end{eqnarray*}
(i) Assume $\beta <0$, then $m=0$ and $l\geq 2$. Since $l''=l-s\geq 1$ holds for all $l\geq 2$, it must be that $s=1$. Thus
$$\frac{l-s+1}{l+1}=1-\frac{s}{l+1}=1-\frac{1}{l+1}\leq 1\ \ for\ \ all\ \ l\geq 2.$$
Let $\rho =1$. Then
$$\int_{0}^{\infty}\frac{r^{l''+1}}{e^{r}}dr\leq \Delta_{0}\rho^{l+1}(l+1)!\ \ if\ \ \int_{0}^{\infty}\frac{r^{l''}}{e^{r}}dr\leq \Delta_{0}\rho^{l}l!.$$
The smallest $l''$ is $l_{0}''=l_{0}-s=2-s$. Now,
\begin{eqnarray*}
  \int_{0}^{\infty}\frac{r^{l_{0}''}}{e^{r}}dr=l_{0}''!=(2-s)!=1! & = & 1\\
                                                                   & = & \frac{1}{2}l!\ \ where\ \ l=2\\
                                                                   & = & \frac{1}{2}\rho^{l}l!\\
                                                                   & = & \Delta_{0}\rho^{l}l!\ \ where\ \ \Delta_{0}=\frac{1}{2}.
\end{eqnarray*} 
It follows that $\int_{0}^{\infty}\frac{r^{l'}}{e^{r}}dr\leq \Delta_{0}\rho^{l}l!$ for all $l\geq 2$.

(ii)Assume $\beta>0$. Then $m=\lceil \frac{\beta}{2}\rceil $ and $l\geq2m+2$. This gives
$$\frac{l-s+1}{l+1}=1-\frac{s}{l+1}\leq 1.$$
Let $\rho =1$. Then $\int_{0}^{\infty}\frac{r^{l''+1}}{e^{r}}dr\leq \Delta_{0}\rho^{l+1}(l+1)!$ if $\int_{0}^{\infty}\frac{r^{l''}}{e^{r}}dr\leq \Delta_{0}\rho^{l}l!$. The smallest $l$ is $l_{0}=2m+2$. Hence the smallest $l''$ is $l_{0}''=l_{0}-s=2m+2-s$. Now,
\begin{eqnarray*}
  \int_{0}^{\infty}\frac{r^{l_{0}''}}{e^{r}}dr & = & l_{0}''!=(2m+2-s)=(l_{0}-s)!\\
                                               & = & \frac{1}{l_{0}(l_{0}-1)\cdots (l_{0}-s+1)}\cdot (l_{0})!\\
                                               & = & \Delta_{0}\rho^{l_{0}}l_{0}!\ where\ \Delta_{0}=\frac{1}{(2m+2)(2m+1)\cdots (2m-s+3)}.    
\end{eqnarray*} 
It follows that $\int_{0}^{\infty}\frac{r^{l'}}{e^{r}}dr\leq \Delta_{0}\rho^{l}l!$ for all $l\geq 2m+2$.\\
\\
{\bf Case3.} Assume $l''=l$. Then
$$\int_{0}^{\infty}\frac{r^{l'}}{e^{r}}dr \leq \int_{0}^{\infty}\frac{r^{l''}}{e^{r}}dr=l!\ \ and\ \ \int_{0}^{\infty}\frac{r^{l'+1}}{e^{r}}dr\leq (l+1)!.$$
Let $\rho=1$. Then
$$\int_{0}^{\infty}\frac{r^{l'}}{e^{r}}dr\leq \Delta_{0}\rho^{l}l!$$
for all $l$ where $\Delta_{0}=1$.
The lemma follows now immediately from the three cases.\hspace{6cm} $\sharp$\\
\\
{\bf Remark.} In the preceding lemma the constant $\rho$ plays an important role in our error estimate. Moreover, the constant $\Delta_{0}$ also means something in the error bound. Therefore we should make a thorough clarification for the two constants. It can be shown easily that
\newcounter{bean}
\begin{list}%
{(\alph{bean})}{\usecounter{bean}}
 \item $l''>l$\ \ if\ \ and\ \ only\ \ if\ \ $n-\beta>3$; 
 \item $l''<l$\ \ if\ \ and\ \ only\ \ if\ \ $n-\beta\leq 1$;
 \item $l''=l$ if\ \ and\ \ only\ \ if\ \ $1<n-\beta\leq 3$
\end{list}
where $l''$ and $l$ are as in the proof of the lemma.

The two constants can then be determined as follows.
\newcounter{book}
\begin{list}%
{(\alph{book})}{\usecounter{book}}
 \item $n-\beta >3$. Let $s=\lceil \frac{n-\beta-3}{2}\rceil$. Then
 \newcounter{car}
 \begin{list}%
 {(\roman{car})}{\usecounter{car}}
 \item if $\beta<0$, then $\rho=\frac{3+s}{3}$ and $\Delta_{0}=\frac{(2+s)(1+s)\cdots 3}{\rho^{2}}$;
 \item if $\beta>0$, then $\rho=1+\frac{s}{2\lceil \frac{\beta}{2}\rceil +3}$ and $\Delta_{0}=\frac{(2m+2+s)(2m+1+s)\cdots(2m+3)}{\rho^{2m+2}}$ where $m=\lceil \frac{\beta}{2}\rceil $.
 \end{list} 
 \item $n-\beta\leq 1$. Let $s=-\lceil \frac{n-\beta-3}{2}\rceil$. Then
 \newcounter{sky}
 \begin{list}%
 {(\roman{sky})}{\usecounter{sky}}
 \item if $\beta<0,$ then $\rho=1$ and $\Delta_{0}=\frac{1}{2}$;
 \item if $\beta>0,$ then $\rho=1$ and $\Delta_{0}=\frac{1}{(2m+2)(2m+1)\cdots (2m-s+3)}$; where $m=\lceil \frac{\beta}{2}\rceil$.
 \end{list}
 \item $1<n-\beta\leq 3$. Then $\rho=1$ and $\Delta_{0}=1$.
\end{list}
Before introducing our main result, we need the following lemma.
\begin{lem}
  For any positive integer $l$,
$$\frac{\sqrt{(2l)!}}{l!}\leq 2^{l}.$$
\end{lem}
{\bf Proof}. This inequality holds for $l=1$ obviously. We proceed by induction.
\begin{eqnarray*}
  \frac{\sqrt{[2(l+1)]!}}{(l+1)!} & = & \frac{\sqrt{(2l+2)!}}{l!(l+1)}=\frac{\sqrt{(2l)!}}{l!}\cdot \frac{\sqrt{(2l+2)(2l+1)}}{l+1}\\
                                  & \leq & \frac{\sqrt{(2l)!}}{l!}\cdot \frac{\sqrt{(2l+2)^{2}}}{l+1}\leq 2^{l}\cdot \frac{(2l+2)}{l+1}=2^{l+1}. 
\end{eqnarray*}
Our lemma thus follows. \hspace{11cm} $\sharp $

Now we have come to our main result.
\begin{thm}
  Let $h$ be as in (4) and $\rho$ be as in Lemma2.1 and the remark following its proof. For any positive real number $b_{0}$, let $C=\max \{ \frac{2}{3b_{0}}, 8\rho \} $ and $\delta_{0}=\frac{1}{3C}$. If $f\in {\cal C}_{h,m}$ and $s$ is the $h$ spline that interpolates $f$ on a subset $X$ of $R^{n}$, as discussed in the beginning of section1, then
\begin{equation}
  |f(x)-s(x)|\leq 2^{\frac{n+\beta-7}{4}}\cdot \pi^{\frac{n-1}{4}}\cdot \sqrt{n\alpha_{n}}\cdot c^{\frac{\beta}{2}}c^{-l}\sqrt{\Delta_{0}}\sqrt{3C}\sqrt{\delta} (\lambda')^{\frac{1}{\delta}}\| f\| _{h}
\end{equation}
holds for all $x\in \Omega$ and $0<\delta <\delta_{0}$ where $0<\lambda' <1$ is a fixed constant and $\Omega$ is any subset of $R^{n}$ satisfying the condition that for any $x\in \Omega$, and any number $r$ with $\frac{1}{3C}\leq r\leq \frac{2}{3C}$(Note that $\frac{2}{3C}\leq b_{0}$), there is a simplex $Q$, $x\in Q\subseteq \Omega$, with $diamQ=r$, such that $Q$ contains equally spaced centers of degree $l$, $\frac{1}{3C\delta}\leq l\leq \frac{2}{3C\delta}$, the centers being contained in $X$. In particular, if $x$ is fixed, the requirement for $\Omega$ is only the existence of a simplex $Q$, $x\in Q \subseteq \Omega$, which satisfies the aforementioned properties, and $X$ can be chosen to consist only of the equally spaced centers in $Q$.

The constants $\beta,\ \alpha_{n}$ and $\Delta_{0}$ are as in (4) and Lemma2.1 and the remark following its proof, and the crucial constant $\lambda'$ is given by 
$$\lambda'=\left( \frac{2}{3}\right) ^{\frac{1}{3C}}$$
and only in some cases mildly depends on the dimension $n$. The $h$-norm of $f$, $\|f\| _{h}$ is defined as in the statement following (7).
\end{thm}
{\bf Proof}. Let $\rho$ be as in Lemma2.1 and $b_{0}$ be arbitrary positive real number. Fix $\delta_{0}=\frac{1}{3C}$. For any $0<\delta\leq \delta_{0}$, it's easily seen that $0<3C\delta \leq 1$. There exists a positive integer $l$ such that $1\leq 3l\delta C\leq 2$. From this we get $\frac{1}{3C\delta}\leq l\leq \frac{2}{3C\delta}$ and $\frac{1}{3C}\leq l\delta \leq \frac{2}{3C}\leq b_{0}$.

For any $x\in \Omega$, let $Q$ be an $n$-simplex containing $x$ such that $Q\subseteq \Omega$ and has diameter $diamQ=l\delta$. Then Theorem4.2 of \cite{MN2} implies that 
\begin{equation}
  |f(x)-s(x)|\leq c_{l}\| f\| _{h}\int _{R^{n}}|y-x|^{l}d|\sigma |(y)
\end{equation}
whenever $l>m$, where $\sigma$ is any measure supported on $X$ such that
\begin{equation}
  \int_{R^{n}}p(y)d\sigma(y)=p(x)
\end{equation}
for all polynomials $p$ in $P_{l-1}^{n}$. Here
$$c_{l}=\left\{ \int_{R^{n}}\frac{|\xi|^{2l}}{(l!)^{2}}d\mu(\xi)\right\} ^{1/2}$$
whenever $l>m$. By (8), for all $2l\geq 2m+2$,
\begin{eqnarray}
  c_{l} & = & \left\{ \int_{R^{n}}\frac{|\xi|^{2l}}{(l!)^{2}}d\mu(\xi)\right\} ^{1/2} \nonumber \\
        & \leq & \frac{1}{l!}2^{\frac{n+\beta+1}{4}}\pi^{\frac{n+1}{4}}\sqrt{n\alpha_{n}}c^{\frac{\beta-2l}{2}}\sqrt{\Delta_{0}}\rho^{l}\sqrt{(2l)!} \nonumber \\
        & \leq & 2^{\frac{n+\beta+1}{4}}\pi^{\frac{n+1}{4}}\sqrt{n\alpha_{n}}c^{\frac{\beta}{2}}c^{-l}\sqrt{\Delta_{0}}(2\rho)^{l}    
\end{eqnarray}
due to Lemma2.2.

In order to obtain the bound on $|f(x)-s(x)|$ as mentioned in the theorem, one only needs to find a suitable bound for
$$I=c_{l}\int_{R^{n}}|y-x|^{l}d|\sigma|(y).$$
This can be attained by choosing the measure $\sigma$ in an appropriate way, as the following deduction shows.

Let $Y=\{ x_{1},\cdots ,x_{N}\} \subseteq X$ be the set of equally spaced centers of degree $l-1$ in $Q$, where $N=dimP_{l-1}^{n}$. By Lemma1.2, there is a measure $\sigma$ supported on $Y$ such that 
$$\int p(y)d\sigma(y)=p(x)$$
for all $p$ in $P_{l-1}^{n}$, and 
$$\int_{Q}d|\sigma|(y)\leq \left( \begin{array}{c}
                                    2l-3 \\ l-1
                                  \end{array} \right) $$  
Now,
\begin{eqnarray*}
  \left( \begin{array}{c}
           2l-3 \\ l-1
         \end{array} \right) = \left( \begin{array}{c}
                                        2(l-1)-1 \\ l-1   
                                      \end{array} \right) & \leq & \left( \begin{array}{c}
                                                                            2(l-1) \\ l-1 
                                                                          \end{array} \right) \\
                                                          & \sim & \frac{1}{\sqrt{\pi}}\cdot \frac{1}{\sqrt{l-1}}\cdot 4^{(l-1)}\ by\ Stirling's\ Formula \\
                                                          & = & \frac{1}{4\sqrt{\pi}}\cdot \frac{1}{\sqrt{l-1}}\cdot 4^{l}.  
\end{eqnarray*}
Therefore, for such $\sigma$, we have
\begin{eqnarray*}
  c_{l}\int_{R^{n}}|y-x|^{l}d|\sigma|(y) & \leq & 2^{\frac{n+\beta +1}{4}}\pi ^{\frac{n+1}{4}}\sqrt{n\alpha_{n}}c^{\frac{\beta}{2}}c^{-l}\sqrt{\Delta_{0}}(2\rho)^{l}(l\delta)^{l}\left( \begin{array}{c}
                                                                                             2l-3 \\l-1
                                                                                           \end{array} \right) \\
                                         & \leq & 2^{\frac{n+\beta+1}{4}}\pi^{\frac{n+1}{4}}\sqrt{n\alpha_{n}}c^{\frac{\beta}{2}}c^{-l}\sqrt{\Delta_{0}}(2\rho)^{l}(l\delta)^{l}\frac{1}{4\sqrt{\pi}}\frac{1}{\sqrt{l-1}}4^{l} \\
                                         & = & 2^{\frac{n+\beta+1}{4}}\pi^{\frac{n+1}{4}}\sqrt{n\alpha_{n}}c^{\frac{\beta}{2}}c^{-l}\sqrt{\Delta_{0}}(8\rho)^{l}(l\delta)^{l}\frac{1}{4\sqrt{\pi}}\frac{1}{\sqrt{l-1}}\\
                                         & \leq & 2^{\frac{n+\beta+1}{4}}\pi^{\frac{n+1}{4}}\sqrt{n\alpha_{n}}c^{\frac{\beta}{2}}c^{-l}\sqrt{\Delta_{0}}\frac{1}{4\sqrt{\pi(l-1)}}(Cl\delta)^{l}\ where\ C=max\left\{ \frac{2}{3b_{0}},\ 8\rho\right\} \\ 
                                         & \leq & 2^{\frac{n+\beta+1}{4}}\pi^{\frac{n+1}{4}}\sqrt{n\alpha_{n}}c^{\frac{\beta}{2}}c^{-l}\sqrt{\Delta_{0}}\frac{1}{4\sqrt{\pi(l-1)}}\left( \frac{2}{3}\right) ^{l}\ \ since\ \ 1\leq 3l\delta C\leq 2 \\ 
                                         & \leq & 2^{\frac{n+\beta+1}{4}}\pi^{\frac{n+1}{4}}\sqrt{n\alpha_{n}}c^{\frac{\beta}{2}}c^{-l}\sqrt{\Delta_{0}}\frac{1}{4\sqrt{\pi(l-1)}}\left[ \left( \frac{2}{3}\right) ^{\frac{1}{3C}}\right] ^{\frac{1}{\delta}}.       
\end{eqnarray*}
Note that $l\rightarrow \infty$ as $\delta\rightarrow 0$. Thus $\sqrt{l-1}\sim \sqrt{l}$ whenever $\delta$ is small. Moreover, $\frac{1}{3C\delta}\leq l\leq\frac{2}{3C\delta}$ implies $\frac{3C\delta}{2}\leq\frac{1}{l}\leq 3C\delta$ and $\sqrt{\frac{3}{2}C\delta}\leq \frac{1}{\sqrt{l}}\leq \sqrt{3C\delta}$. Thus, for $\delta$ small enough, we get 
\begin{eqnarray*}
|f(x)-s(x)| & \leq & 2^{\frac{n+\beta+1}{4}}\pi^{\frac{n+1}{4}}\sqrt{n\alpha_{n}}c^{\frac{\beta}{2}}c^{-l}\sqrt{\Delta_{0}}\frac{1}{4\sqrt{\pi}}\sqrt{3C}\sqrt{\delta}(\lambda')^{\frac{1}{\delta}}\|f\| _{h}\\
            & = & 2^{\frac{n+\beta-7}{4}}\pi^{\frac{n-1}{4}}\sqrt{n\alpha_{n}}c^{\frac{\beta}{2}}c^{-l}\sqrt{\Delta_{0}}\sqrt{3C}\sqrt{\delta}(\lambda')^{\frac{1}{\delta}}\| f\| _{h}
\end{eqnarray*}
where $\lambda'=\left( \frac{2}{3}\right) ^{\frac{1}{3C}}$. \hspace{12cm}\ \ \ \ \     $\sharp$\\
\\
{\bf Remark}. In this theorem although we don't explicitly state that the number $\delta$ represents the well-known fill-distance $d(Q,Y)$, it is equivalent to the fill-distance in spirit. This can be seen from the inequality $\frac{1}{3C}\leq l\delta\leq \frac{2}{3C}$ where $l$ denotes the degree of the equally spaced centers $Y$ in $Q$. Note that $\delta\rightarrow 0$ iff $l\rightarrow \infty$ and $l\rightarrow \infty$ iff $d(Q,Y)\rightarrow 0$. In this paper we avoid using the term fill-distance since the data points are not purely scattered. However the equally spaced centers in a simplex are theoretically very easy to handle and in application impose no trouble at all on programming. Moreover, in the proof of the theorem, $\delta$ is usually small. Therefore the gap between $\sqrt{l-1}$ and $\sqrt{l}$ can be ignored. This is why we adopt ``$\leq$'' in (9).

The second thing is about the function space. Although the space ${\cal C}_{h,m}$ adopted by us is defined by Madych and Nelson in \cite{MN1} and \cite{MN2} and is not the same as Wu and Schaback's function space \cite{Wu}, there is a nice unification theory in \cite{We} which is achieved by putting mild restrictions on Wu and Schaback's function space. We call it {\bf native space}.
\section{Comparison}
The original exponential-type error bound for multiquadric interpolation (See \cite{MN3} and \cite{Lu3}.) is of the form
\begin{equation}
  |f(x)-s(x)|\leq C_{0}\lambda^{\frac{1}{d}}\| f\| _{h},\ d\rightarrow 0
\end{equation}
, where $C_{0}$ is about the same as the long string of coefficients in (9), $\lambda$ is a constant satisfying $0<\lambda <1$, and $d$ is the fill-distance equivalent to $\delta$ of (9). The definition of $\lambda$ is 
$$\lambda =\left( \frac{2}{3} \right) ^{\frac{1}{3C\gamma_{n}}}$$
where
$$C=\max \left\{ 2\rho'\sqrt{n}e^{2n\gamma_{n}},\ \frac{2}{3b_{0}}\right\} ,\ \rho'=\frac{\rho}{c}$$
with $\rho$ and $c$ the same as this paper, $b_{0}$ is the side length of a cube, equivalent to the $b_{0}$ of Theorem2.3 of this paper, and $\gamma_{n}$ is defined recurvely by $\gamma_{1}=2,\ \gamma_{n}=2n(1+\gamma_{n-1})$ if $n>1$. The first few $\gamma_{n}$ are
$$\gamma_{1}=2,\ \gamma_{2}=12,\ \gamma_{3}=78,\ \gamma_{4}=632,\ \gamma_{5}=6330, \cdots $$. The fast growth of $\gamma_{n}$ as $n\rightarrow \infty$ forces $e^{2n\gamma_{n}}$ to grow extremely rapidly, driving $\lambda$ to 1 and making the error bound (13) meaningless.

In our new approach, the crucial constant $\lambda'$ in (9) is defined by 
$$\lambda'=\left( \frac{2}{3} \right) ^{\frac{1}{3C}}$$
where
$$C=\max \left\{ \frac{2}{3b_{0}},\ 8\rho \right\}$$
. Note that $\rho$ is only mildly dependent of dimension $n$. This can be seen in Lemma2.1 and the remark following its proof. A lot of time $\rho$, and hence $\lambda'$, is independent of $n$. For example, when $\beta=1$ and $n=1,2,3$ or $4$, $\lambda'$ is completely independent of $n$. Fortunately, these are exactly the most interesting and important cases. For higher dimensions, $\lambda'$ can remain unchanged by increasing $\beta$ in (4) to keep $n-\beta \leq 3$.

Besides this, $\sqrt{\delta}$ in (9) also contributes to the convergence rate of the error bound as $\delta\rightarrow 0$.

\end{document}